\documentclass[12pt]{amsart}

\usepackage{graphicx}
\usepackage{caption}
\usepackage{subcaption}
\usepackage{color}
\usepackage{amsmath,amsfonts}
\usepackage{algorithm}
\usepackage{algorithmic}

\theoremstyle{definition}

\theoremstyle{remark}

\begin{document}


\title[AMG for quadratic FEM equations]{An algebraic multigrid method for quadratic finite element equations of elliptic and saddle point systems in 3D}



\author{Huidong Yang}
\address{Johann Radon Institute for Computational and Applied Mathematics (RICAM), Austrian Academy of Sciences, Altenberger Strasse 69, A-4040 Linz, Austria}
\email{huidong.yang@oeaw.ac.at}



\begin{abstract}
  In this work, we propose a robust and easily implemented algebraic multigrid 
  method as a stand-alone solver or a preconditioner in Krylov subspace methods 
  for solving either symmetric and positive definite or saddle point 
  linear systems of equations arising from the finite element discretization of the 
  vector Laplacian problem, linear elasticity problem in 
  pure displacement and mixed displacement-pressure form, 
  and Stokes problem in mixed velocity-pressure form in 3D, 
  respectively. 
  We use hierarchical quadratic basis functions to construct the finite element 
  spaces. A new heuristic algebraic coarsening strategy is introduced for 
  construction of the hierarchical coarse system matrices. 
  We focus on numerical study of the mesh-independence robustness 
  of the algebraic multigrid and the algebraic multigrid 
  preconditioned Krylov subspace methods. 
\end{abstract}


\maketitle









\section{Introduction}
Compared to the geometrical multigrid (GMG) method (see, .e.g., \cite{HB03}), 
the algebraic multigrid (AMG) method (see, e.g., \cite{KS01:00, JWRKS87}) is a purely 
matrix-based approach, that does not rely on any underlying mesh 
hierarchy; see, e.g., \cite{UH:02} for the development 
from GMG to AMG methods. 
Concerning comparison of different types of AMG methods  we refer to, e.g., 
\cite{Stuben2001281} for a review and related references. 
In contrast to coarsening based on the strongly connected matrix entries in the 
classical AMG method, an 
AMG method (among others) with special coarsening and 
interpolation strategies was introduced in \cite{FK98:00}, that is based on 
graph connectivity of the matrix only and leads to 
fast construction of matrices on coarse levels. 
An AMG method that is based on the matrix graph information only 
was also studied early in \cite{DB95}. The further development 
of such an AMG method \cite{FK98:00} in different applications have been 
reported in, e.g., \cite{NME:NME839,WM04:00,WM06:00,Langer2005406,Langer201547,Yang20115367}. In this work, we focus on the development of such an AMG method 
for both elliptic and saddle point systems of equations arising from the quadratic finite element 
discretization for the three dimensional (3D) 
vector Laplacian problem, linear elasticity problem in pure displacement and 
mixed displacement-pressure form, and the Stokes problem 
in mixed velocity-pressure form. This requires new coarsening strategies to 
construct the hierarchy of matrices on coarse levels for both the elliptic and 
saddle point systems, that are to be developed in this work. We notice, that 
different AMG methods towards higher-order finite element 
equations for second order elliptic problems were also studied using different approaches 
in, e.g., \cite{XUJ06, NOT14}. The main focus of this work is the numerical 
study of the robustness and efficiency of the designed AMG method as a stand-alone solver 
or a preconditioner in Krylov subspace methods for solving the elliptic and saddle point 
systems. 

The remainder of this paper is organized in the following way. 
In Section \ref{sec:prelim}, we describe the model problems, 
their finite element discretizations and the arising linear systems of equations. The 
algebraic multigrid method using a new heuristic coarsening strategy is 
prescribed in Section \ref{sec:amg}. In Section \ref{sec:numres}, we present 
numerical results of the AMG method applied to discrete model problems. Finally, some 
conclusions are drawn in Section \ref{sec:con}.

\section{Preliminaries}\label{sec:prelim}
\subsection{The model problems}
Let $\Omega\subset\mathbb{R}^3$ be a simply connected and bounded domain 
with two boundaries $\Gamma_N$ and $\Gamma_D$ such that 
$\bar{\Gamma}_D\cup\bar{\Gamma}_N=\partial\Omega$ and 
$\Gamma_N\cap\Gamma_D=\emptyset$. We consider the 3D 
vector Laplacian problem, the linear elasticity problem in pure 
displacement and mixed displacement-pressure forms, and the Stokes problem in 
mixed velocity-pressure form, that are formulated in the following:

For the vector Laplacian problem: Find the potential 
$u:\bar{\Omega} \mapsto \mathbb{R}^3$ such that 
\begin{equation}\label{eq:veclap}
  -\Delta u = 0 {\textup{ in }}\Omega 
\end{equation}
with the boundary conditions $u=g_D$ on $\Gamma_D$ 
and $\frac{\partial u}{\partial n}=g_N$ on 
$\Gamma_N$, where $n$ denotes the outward normal vector on $\Gamma_N$. 

For the linear elasticity problem in pure displacement form: Find the 
displacement $u:\bar{\Omega}\mapsto \mathbb{R}^3$ such that 
\begin{equation}\label{eq:elasdis}
   -\nabla\cdot\sigma(u) = 0 {\textup{ in }}\Omega 
\end{equation}
with the boundary conditions $u=g_D$ on $\Gamma_D$ and $\sigma(u)n=g_N$ on 
$\Gamma_N$. In particular, we use the linear Saint Venant-Krichoff elasticity model. 
The Cauchy stress tensor and  the 
infinitesimal strain tensor are defined by  
$\sigma(u) = 2\mu \varepsilon(u) + \lambda \text{div}(u) I$ 
and $\varepsilon(u) = (\nabla u + \nabla u^T)/2$, respectively, 
with Lam\'{e} constants $\lambda$ and $\mu$.

For the linear elasticity problem in mixed displacement-pressure form: Find the 
displacement $u:\bar{\Omega}\mapsto \mathbb{R}^3$ and pressure 
$p:\bar{\Omega}\mapsto \mathbb{R}$ 
such that 
\begin{equation}\label{eq:elasmix}
  \begin{aligned}
    -\nabla\cdot(2\mu\varepsilon(u)) + \nabla p = 0 & {\textup{ in }}\Omega \\
    -\nabla\cdot u -\frac{1}{\lambda} p =  0 & {\textup{ in }}\Omega \\
   \end{aligned}
\end{equation}
with the boundary conditions $u=g_D$ on $\Gamma_D$ and $(2\mu\varepsilon(u)-pI)n=g_N$ on 
$\Gamma_N$. It is easy to see that, in this classical mixed displacement-pressure 
form, the displacement and pressure are associated by the relation 
$p= - \lambda \nabla\cdot u $; see, e.g., \cite{DB07:00}. 

For the Stokes problem in mixed velocity-pressure form: Find the 
velocity $u:\bar{\Omega} \mapsto \mathbb{R}^3$ and pressure 
$p:\bar{\Omega} \mapsto \mathbb{R}$ 
such that 
\begin{equation}\label{eq:stokesmix}
  \begin{aligned}
    -\nabla\cdot(2\mu\varepsilon(u)) + \nabla p = 0 & {\textup{ in }}\Omega \\
    -\nabla\cdot u =  0 & {\textup{ in }}\Omega \\
   \end{aligned}
\end{equation}
with the boundary conditions $u=g_D$ on $\Gamma_D$ 
and $(2\mu\varepsilon(u)-pI)n=g_N$ on 
$\Gamma_N$, where $\mu$ denotes the dynamic viscosity. 

\subsection{The variational formulations}
We search for weak solutions of the above four model problems 
(\ref{eq:veclap})-(\ref{eq:stokesmix}) in proper spaces. For this, 
let $H^1(\Omega)$ and $Q=L^2(\Omega)$ denote the standard 
Sobolev and Lebesgue spaces on $\Omega$; see \cite{AF03:00}. 
With $V={H^1(\Omega)}^3$, we define 
the spaces 
$V_g=\{u\in V: u|_{\Gamma_D}= g_D \}$ 
for the potential and displacement (velocity) functions. We also define the 
homogenized space $V_0=\{u\in V : u|_{\Gamma_D}=0 \}$. 
In addition, we assume the given data 
$g_D\in {H^{1/2}(\Gamma_D)}^3$, where ${H^{1/2}(\Gamma_D)}^3$ denotes the 
trace space, i.e., ${H^{1/2}(\Gamma_D)}^3=\{v|_{\Gamma_D} : v\in H^1(\Omega)^3\}$. 
We also assume the given data $g_N\in L^2(\Gamma_N)^3$. By standard techniques, the following 
variational formulations are obtained. 

The variational formulation for the vector Laplacian problem (\ref{eq:veclap}) and 
the linear elasticity problem (\ref{eq:elasdis}) 
in pure displacement form reads (after homogenization): Find $u\in V_0$ such that
\begin{equation}\label{eq:weakveclapform}
  a( u , v ) = \langle F , v \rangle
\end{equation}
for all $v\in V_0$, with the bilinear form 
$ a( u , v ):=\int_{\Omega} \nabla u:\nabla v dx$ for the vector Laplacian problem and 
$a( u , v ) := \int_{\Omega}[ 2\mu\varepsilon(u):\varepsilon(v) + 
\lambda\nabla\cdot u \nabla\cdot v]dx$ for the linear elasticity problem, respectively, 
and the linear form 
$\langle F , v \rangle := \int_{\Gamma_N} g_N\cdot v dx-a( g_D , v )$, accordingly. 


The variational formulation for the linear elasticity problem (\ref{eq:elasmix}) 
in mixed displacement-pressure form and the Stokes problem (\ref{eq:stokesmix}) in 
mixed velocity-pressure form reads (after homogenization):
Find $u\in V_0$ and $p\in Q$ such that
\begin{equation}\label{eq:weakelasmixform}
  \begin{aligned}
    & a(u , v ) + b( v, p ) = \langle F , v \rangle ,\\
    & b(u, q ) -c(p , q )  =\langle G , q \rangle
  \end{aligned}
\end{equation}
for all $v\in V_0$ and $q\in Q$, where the bilinear and linear forms are given by 
$a(u , v) = 2\mu\int_{\Omega} \varepsilon(u):\varepsilon(v) dx$, 
$b( v, q ) = -\int_{\Omega} q \nabla\cdot v dx $, 
$\langle F , v \rangle = \int_{\Gamma_N} g_N\cdot v dx - a(g_D , v)$, 
$\langle G , q \rangle = 0 $ , respectively, and 
$c(p, q) = \frac{1}{\lambda}  \int_{\Omega}p q dx$ and $c(p, q)=0$ for the linear elasticity 
and Stokes problem, respectively. 

\subsection{The finite element discretization}
The spatial discretization is done by the Galerkin finite element 
method with  a hierarchical quadratic polynomial basis functions. 
Let ${\mathcal T}_h$ be the admissible subdivision of the domain $\Omega$ into 
tetrahedra. The four linear 
basis functions on each tetrahedron $T\in {\mathcal T}_h$ 
are nothing but standard $P_1$ hat functions in 3D, i.e., $\phi_i:=\lambda_i$, 
$i=1 , ... , 4$, where $\lambda_i$ are the 
barycentric coordinates of $T$. The six 
quadratic basis functions are then defined as 
\begin{equation*}
  \begin{aligned}
    &\phi_5= 4\lambda_1\lambda_2,\; 
    \phi_6= 4\lambda_2\lambda_3,\; 
    \phi_7= 4\lambda_3\lambda_4,\;  \\
    &\phi_8= 4\lambda_1\lambda_3,\; 
    \phi_9= 4\lambda_1\lambda_4,\; 
    \phi_{10}= 4\lambda_2\lambda_4,\; 
      \end{aligned}
\end{equation*}
that construct hierarchical quadratic polynomial basis functions. 

Let $V_L:= \{v\in C_0(\bar{\Omega}) : v_{T}\in\Phi_T, \forall T\in {\mathcal T}_h\}$ 
be the subspace of continuous piecewise linear hat functions with zero 
traces on $\Gamma_D$ and 
$V_Q:= \{v\in C_0(\bar{\Omega}) : v_{T}\in\Psi_T, \forall T\in {\mathcal T}_h\}$ 
the subspace of continuous piecewise quadratic functions with zero traces 
on $\Gamma_D$, where $\Phi_T=\text{span}\{\phi_{i=1,...,4} : T\in {\mathcal T}_h\}$ and 
$\Psi_T=\text{span}\{\phi_{i=5,...,10} : T\in {\mathcal T}_h \}$.  
Let $Q_L:= \{v\in C(\bar{\Omega}) : v_{T}\in\Phi_T, \forall T\in {\mathcal T}_h\}$ 
be the subspace of continuous piecewise linear hat functions. 
It can be shown that the global degrees of freedom (DOF) of $V_L$,  $V_Q$, $Q_L$ 
are the number of vertices ($l$) and edges ($m$) (excluding the 
vertices and edges on $\Gamma_D$), and the number of all vertices ($n$), 
respectively. The global basis 
functions $\varphi_i\in V_L, Q_L$ and $\psi_i\in V_Q$ 
can be constructed from the local ones. 
The function space for one component of the potential or displacement is 
defined as $V_h=V_L\oplus V_Q\subset V_0$, a linear subspace complemented 
by a quadratic subspace. 
The function space for the pressure is defined as $Q_h=Q_L\subset Q $. 

Using Galerkin's principle the discrete 
elliptic variational formulation for the vector Laplacian and linear elasticity problem 
in pure displacement form read: 
Find $u_h\in V_h$ such that 
\begin{equation}\label{eq:disveclapform}
  a( u_h , v_h ) = \langle F , v_h \rangle
\end{equation}
for all $v_h\in V_h$. 

The discrete mixed variational formulation for the elasticity problem 
in mixed displacement-pressure form and the Stokes problem in mixed 
velocity-pressure form reads: Find $(u_h , p_h)\in V_h\times Q_h$ such that 
\begin{equation}\label{eq:diselasmixform}
  \begin{aligned}
    & a(u_h , v_h ) + b( v_h , p_h ) = \langle F , v_h \rangle ,\\
    & b(u_h, q_h ) -c(p_h , q_h )  =\langle G , q_h \rangle
  \end{aligned}
\end{equation}
for all $v_h\in V_h$ and $q_h\in Q_h$. 

The finite element solutions $u_h\in V_h$ and $p_h\in Q_h$ 
are expressed by the ansatz:
\[
u_h = u_h^l + u_h^q = \displaystyle\sum_{i=1}^{l} u_i^l\varphi_i + 
\displaystyle\sum_{i=1}^{m} u_i^q\psi_i ,\quad
p_h = \displaystyle\sum_{i=1}^{n} p_i^l\varphi_i,
\]
respectively, where $u_i^l, u_i^q \in \mathbb{R}^3$ and $p_i^l\in \mathbb{R}$. 
It is easy to see the finite element 
solution $u_h$ is the sum of the linear and quadratic part, $u_h^l$ 
and $u_h^q$, respectively. 
For the pressure $p_h$, we have a linear approximation. The mixed 
finite element for the elasticity and Stokes problem 
is classical Taylor-Hood element, that fulfills 
the $\inf-\sup$ stability requirement; see, e.g., \cite{BF91:00}.

\subsection{SPD and saddle point linear systems of equations}
Using the finite element discretization (including homogenization), 
we obtain the following symmetric and positive definite (SPD) system of equations for the 
elliptic problem:
\begin{equation}\label{eq:ellipticeqn}
  Au =  \left[ 
    \begin{array}{cc}
      K_{ll} & K_{ql}^T\\
      K_{ql} & K_{qq} \\
    \end{array}
    \right]
  \left[
    \begin{array}{c}
      \underline{u}_l \\ \underline{u}_q
    \end{array}
    \right]
  =
  \left[
    \begin{array}{c}
      \underline{f}_l \\ \underline{f}_q 
    \end{array}
    \right]
  = f,
\end{equation}
where $K_{ll}=(a(\varphi_i, \varphi_j))$, $K_{ql}=a(\varphi_j, \psi_i)$, $K_{qq}=a(\psi_i, \psi_j)$, 
$\underline{u}_l=(u_i^l)$, $\underline{u}_q=(u_i^q)$, 
$\underline{f}_l=(\langle F , \varphi_i\rangle)$ 
and $\underline{f}_q=(\langle F , \psi_i\rangle)$.

For the elasticity problem in mixed displacement-pressure form and the Stokes 
problem in mixed velocity-pressure form, we obtain 
the following symmetric indefinite system of equations:
\begin{equation}\label{eq:mixedeqn}
  \underbrace{\left[ 
     \begin{array}{cc}
      A & B^T\\
      B& -C
     \end{array}
    \right]}_{=:K}
  \left[ 
     \begin{array}{c}
       u \\ p
     \end{array}
    \right]=
  \left[ 
    \begin{array}{cc|c}
      K_{ll} & K_{ql}^T & B^T_{ll}\\
      K_{ql} & K_{qq} & B^T_{lq}\\ \hline
      B_{ll} & B_{lq} & -C_{ll}
    \end{array}
    \right]
  \left[
    \begin{array}{c}
      \underline{u}_l \\ \underline{u}_q \\ \hline \underline{p}_l
    \end{array}
    \right]
  =
  \left[
    \begin{array}{c}
      \underline{f}_l \\ \underline{f}_q \\ \hline \underline{g}_l 
    \end{array}
    \right]
  =
  \left[
    \begin{array}{c}
      f \\ g
    \end{array}
    \right]
\end{equation}
where $K_{ll}=(a(\varphi_i, \varphi_j))$, $K_{ql}=a(\varphi_j, \psi_i)$, 
$K_{qq}=a(\psi_i, \psi_j)$, 
$B_{ll}=b(\varphi_i, \varphi_j)$, $B_{lq}=b(\varphi_i, \psi_j)$, 
$C_{ll}=c(\varphi_i, \varphi_h)$, $\underline{u}_l=(u_i^l)$, $\underline{u}_q=(u_i^q)$, 
$\underline{p}_l=(p_i^l)$, 
$\underline{f}_l=(\langle F , \varphi_i\rangle)$, 
$\underline{f}_q=(\langle F , \psi_i\rangle)$ and 
$\underline{g}_l=(\langle G , \varphi_i \rangle)$. It is obvious 
that for the Stokes problem, $C=0$. 

In the following, we focus on how to solve the above two systems of equations 
(\ref{eq:ellipticeqn}) and (\ref{eq:mixedeqn}) using an algebraic multigrid method.

\section{An algebraic multigrid method}\label{sec:amg}

\subsection{The basic algebraic multigrid iteration}
The basic AMG iteration applied to a general linear system of 
equations $Kx=b$ is given in Algorithm \ref{alg:amgit}, 
with $m_{pre}$ and $m_{post}$ being the 
number of pre- and post-smoothing steps (steps 1-3 and 14-16, respectively). 
By choosing $\nu=1$ and 
$\nu=2$, the iterations in Algorithm \ref{alg:amgit} are called 
V- and W-cycle, respectively. As 
a convention, we use $l=0,..., L$ to indicate the 
algebraic multigrid levels from the finest level $l=0$ to the 
coarsest level $l=L$. 
On the coarsest level $L$, the system is solved by any direct solver (step 6). 
The coarse grid correction step is indicated in steps 4-13. 
The full AMG iterations are realized by repeated application 
of this algorithm. The iteration in this algorithm 
is also combined with the Krylov subspace methods, 
that usually leads to accelerated convergence of 
V-cycle or W-cycle preconditioned methods \cite{Stuben2001281}. 
\begin{algorithm}
   \caption{Basic AMG iteration: AMG($K_l, x_l, b_l$) }\label{alg:amgit} 
   \begin{algorithmic}[1]
     \FOR{$k=1$ to $m_{pre}$} 
     \STATE $x_l^{k+1}={\mathcal S}_l(x_l^k, b_l)$
     \ENDFOR
     \STATE $b_{l+1}=R_l^{l+1}(b_l-K_lx_l)$,
     \IF{l+1=L}
     \STATE Solve $K_Lx_L=b_L$
     \ELSE
     \STATE{
       $x_{l+1}=0$,
       \FOR{$k=1,...,\nu$}
       \STATE $x_{l+1}=$AMG($K_{l+1}, x_{l+1}, b_{l+1}$),
       \ENDFOR
     }
     \ENDIF
     \STATE $x_l=x_l+P_{l+1}^lx_{l+1}$,
     \FOR{$k=1$ to $m_{post}$} 
     \STATE $x_l^{k+1}={\mathcal S}_l(x_l^k, b_l)$,
     \ENDFOR
     \STATE return $x_l$.
   \end{algorithmic}
\end{algorithm}
\subsection{A new heuristic coarsening strategy}
\subsubsection{Case I : The SPD system}
A robust coarsening 
strategy is an important feature of the AMG method, that is 
used to construct the system matrices on coarse levels $l=1, 2, ..., L$. For the 
second order elliptic equations discretized by low order finite element 
or boundary element methods, 
some well known 
graph-based black-box or grey-box type AMG methods have been introduced and applied, 
see, e.g., \cite{DB95, FK98:00,  SR01:00,  NME:NME839, Langer2005406}. 
The general strategy is to split the nodes into the sets of coarse and fine nodes, 
based on the graph connectivity of the system matrix or the constructed auxiliary 
matrix ("virtual'' finite element mesh \cite{SR01:00}). 

\begin{figure}[htbp]
  \centering
  \includegraphics[scale=0.6]{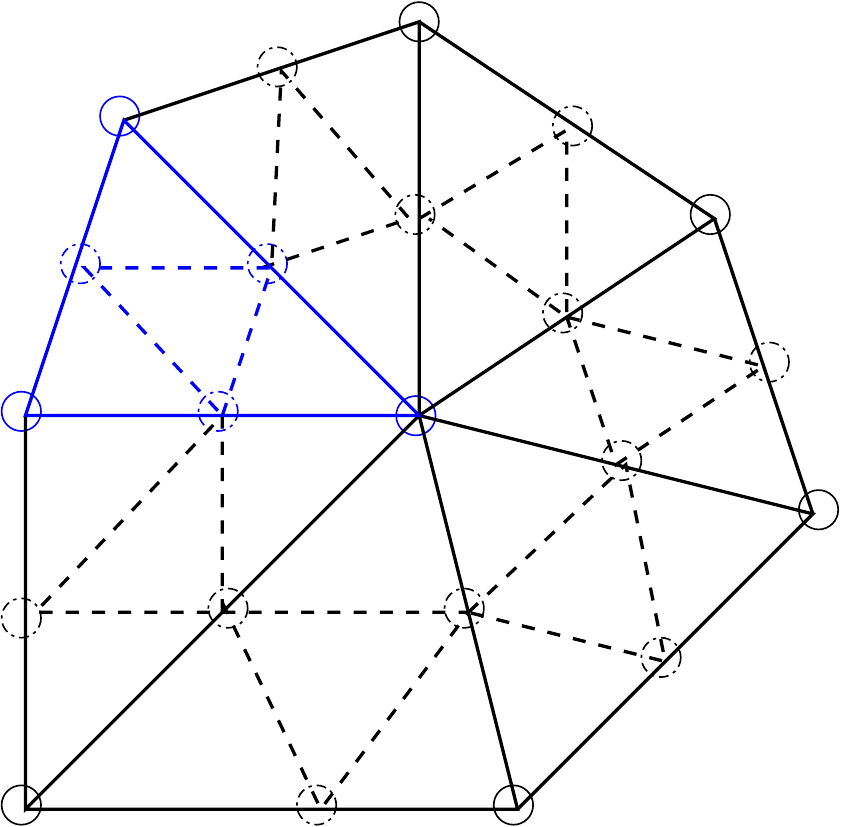}\quad
  \includegraphics[scale=1.2]{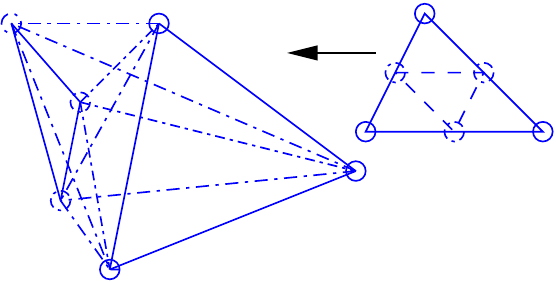}
  \caption{Graph connectivity constructed using the new (left) and the classical (right) 
    strategies: linear DOF (solid line) , quadratic DOF (dashed line) , 
    coupling (dashed dot lines).}\label{fig:amgellipcoarse}
\end{figure}

When applying such a technique to the system matrix in (\ref{eq:ellipticeqn}), we 
obtain a very dense graph connectivity constructed from the stiffness matrix $A$, 
that contains connectivities for the linear DOF, the quadratic DOF and the 
coupling between them. This will lead to a mixture of different oder of DOF 
and may cause additional difficulty to construct the interpolation operators. In fact, from 
our numerical studies, we observe the loss of optimality of the AMG method when 
such a dense graph connectivity is adopted for coarse system matrix construction. 
Therefore, we construct the graph connectivities only for the linear and quadratic DOF, i.e., 
for $K_{ll}$  and $K_{qq}$, respectively. We simply neglect the 
coupling connectivity in $K_{ql}$. By this means, we avoid the mixture of different order of 
DOF on the coarse level. In addition, we are able to construct and control the 
interpolation operators for the linear and quadratic part, respectively. 
A simple comparison of the classical 
and new graph connectivities is illustrated in Fig. \ref{fig:amgellipcoarse}. 

From the left plot in Fig. \ref{fig:amgellipcoarse}, we show the two graph connectivities 
indicated by solid and dashed lines for the linear and quadratic DOF, respectively. 
For a comparison, on the right plot, 
the graph connectivities constructed by the new strategy is reconstructed 
by the classical strategy. For simplicity, we only reconstruct the part indicated by 
the lines with blue color. It is easy to obverse that, the classical one leads to much denser 
graph connectivities than the new one due to the coupling. 

Based on these two graph connectivities, the prolongation matrix $P_{l+1}^l$ 
from the coarse level $l+1$ to the next finer level $l$ is constructed in form of 
\begin{equation}\label{eq:disprolongmat}
P_{l+1}^l =  \left[\begin{array}{cc} I_{l+1}^l & \\ & J_{l+1}^l \\ \end{array} \right],
\end{equation}
where the prolongation matrices, 
$I_{l+1}^l : ({\mathbb R^3})^{n_{l+1}}\rightarrow({\mathbb R}^3)^{n_{l}}$ and 
$J_{l+1}^l : ({\mathbb R^3})^{m_{l+1}}\rightarrow({\mathbb R}^3)^{m_{l}}$ are defined 
for the linear and quadratic DOF, respectively, where $n_{l}$ and $m_l$ 
denote the number of linear and quadratic DOF on level $l$, respectively. 
The restriction matrix from the finer level $l$ to the 
next coarser level $l+1$ is constructed as $(P_{l+1}^l)^T$. 
The system matrix $A_{l+1}$ on the level $l+1$ is constructed by the Galerkin 
projection method:
\[
A_{l+1} = (P_{l+1}^l)^TA_l P_{l+1}^l.
\]
We mention that the two graph connectivities for the linear and 
quadratic DOF are naturally different as illustrated in Fig. \ref{fig:amgellipcoarse}. 
that may require different coarse and fine nodes selection algorithms. However, for 
simplicity, we apply the coarse and fine nodes 
selection and prolongation operator matrix construction 
algorithms developed in \cite{FK98:00} 
for both the linear and quadratic DOF, that show the robustness from the numerical 
studies. 

\subsubsection{Case II : The saddle point system}
A robust coarsening strategy for saddle point problems is, in general, more 
involved than that for elliptic problems mainly due to the $\inf-\sup$ instability issues possibly 
caused 
by standard Galerkin projection method. For saddle point problems arising from the 
low order finite element discretized fluid problem, the stability issue has been 
studied in, e.g., \cite{WM06:00, WM04:00, BM:13}.  In \cite{MW03:00, WM04:00}, 
a so-called {\it 2-shift} coarsening strategy was introduced for the discrete fluid problem 
using the modified Taylor-Hood element ($P_1$iso$P_2-P_1$), that mimics 
the hierarchy of matrices in the geometrical multigrid method. To guarantee 
the $\inf-\sup$ stability for the coarse system is still a research topic. 

\begin{figure}[htbp]
  \centering
  \includegraphics[scale=0.6]{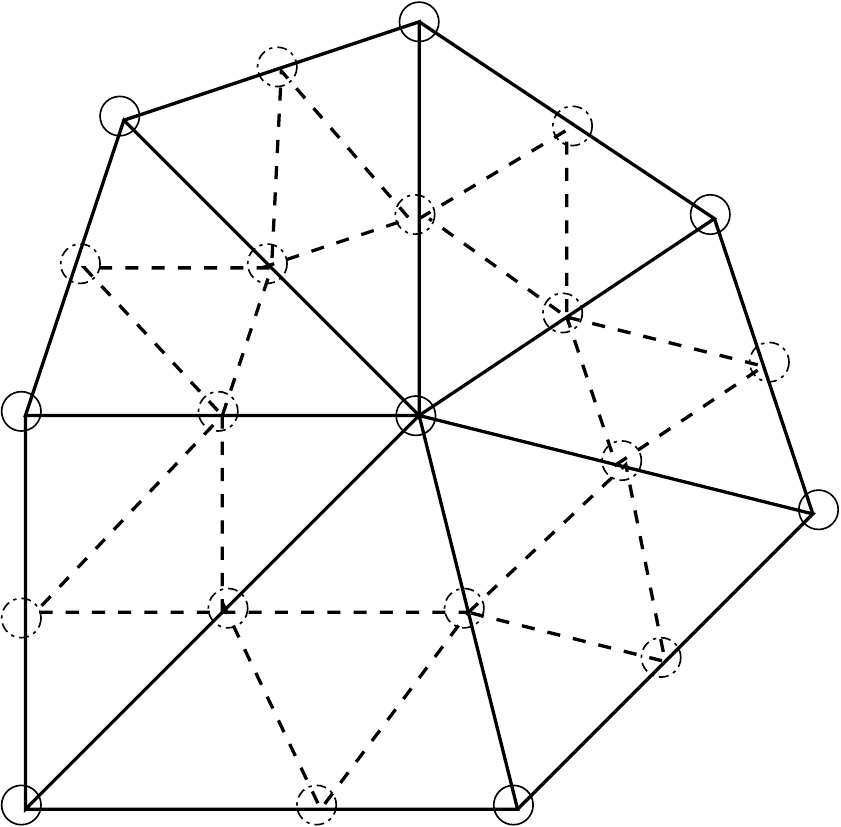}
  \includegraphics[scale=0.6]{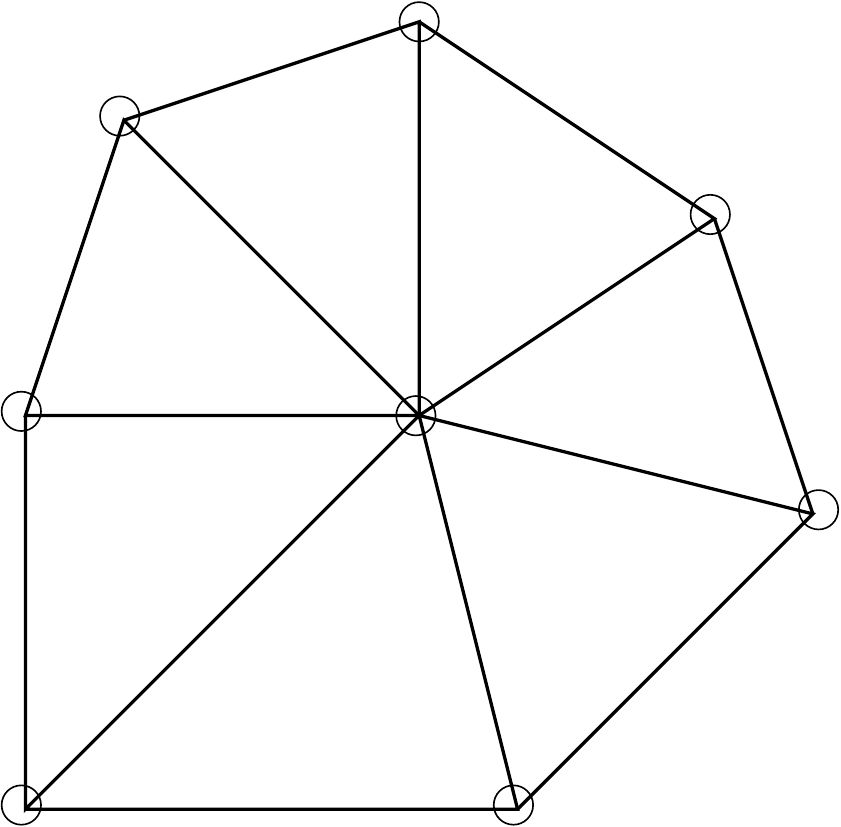}
  \caption{Graph connectivities constructed using the new strategies for 
    velocity (displacement) (left) and pressure (right).}\label{fig:amgmixcoarse}
\end{figure}

We extend the new coarsening strategy described above for the elliptic problem, 
to the saddle point problem, based on the new graph connectivity construction. 
As illustrated in Fig. \ref{fig:amgmixcoarse}, we show the graph connectivities constructed 
for the velocity (displacement) and pressure. For the velocity, we follow the same 
strategy as for the SPD system; see the left plot. 
For the pressure, we have conventional graph 
connectivity for the linear finite element matrix; see the right plot.  

Based on the graph connectivities, the prolongation matrix $P_{l+1}^l$ 
from the coarse level $l+1$ to the next finer level $l$ is constructed in form of 
\begin{equation}\label{eq:mixprolongmat}
P_{l+1}^l =  \left[\begin{array}{ccc} I_{l+1}^l & &\\ & J_{l+1}^l & \\ && H_{l+1}^l\end{array} \right],
\end{equation}
where the prolongation matrices, 
$I_{l+1}^l : ({\mathbb R^3})^{n_{l+1}}\rightarrow({\mathbb R}^3)^{n_{l}}$ and 
$J_{l+1}^l : ({\mathbb R^3})^{m_{l+1}}\rightarrow({\mathbb R}^3)^{m_{l}}$ are defined 
for the linear and quadratic velocity DOF, respectively, where $n_{l}$ and $m_l$ 
denote the number of linear and quadratic velocity DOF on level $l$, respectively, 
$H_{l+1}^l : {\mathbb R}^{k_{l+1}}\rightarrow{\mathbb R}^{k_{l}}$ 
for the pressure DOF, where $k_{l}$ the number of pressure DOF on level $l$. 
The restriction matrix from the finer level $l$ to the 
next coarser level $l+1$ is constructed as $(P_{l+1}^l)^T$. 
The system matrix $K_{l+1}$ on the level $l+1$ is constructed by the Galerkin 
projection method:
\[
K_{l+1} = (P_{l+1}^l)^TK_l P_{l+1}^l.
\]
We admit that the $\inf-\sup$ stability of the coarse system is still open 
by this construction. However, from the numerical studies, 
we observe quite satisfactory results using this new coarsening strategy. 
Nevertheless, we are at least able to obtain efficient multigrid preconditioners by 
using pure Galerkin projection; see comments in, e.g., \cite{KS01:00, TW15} and 
numerical experiments for the Stokes problem using the low order finite element 
discretization in, e.g., \cite{JA08:00}.

\subsection{The smoothing procedure}
To complete the algebraic multigrid algorithm, a smoothing 
procedure is needed. As conventional choices, 
we employ the damped block Jacobi and block Gauss-Seidel smoothers 
for the SPD system (\ref{eq:ellipticeqn}), that 
are widely used in the multigrid methods.  For the saddle point 
system (\ref{eq:mixedeqn}), we have considered the following 
smoothers, that were originally designed and analyzed 
in the GMG method. 
\subsubsection{The multiplicative Vanka smoother}
\begin{figure}[htbp]
  \centering
  \includegraphics[scale=0.8]{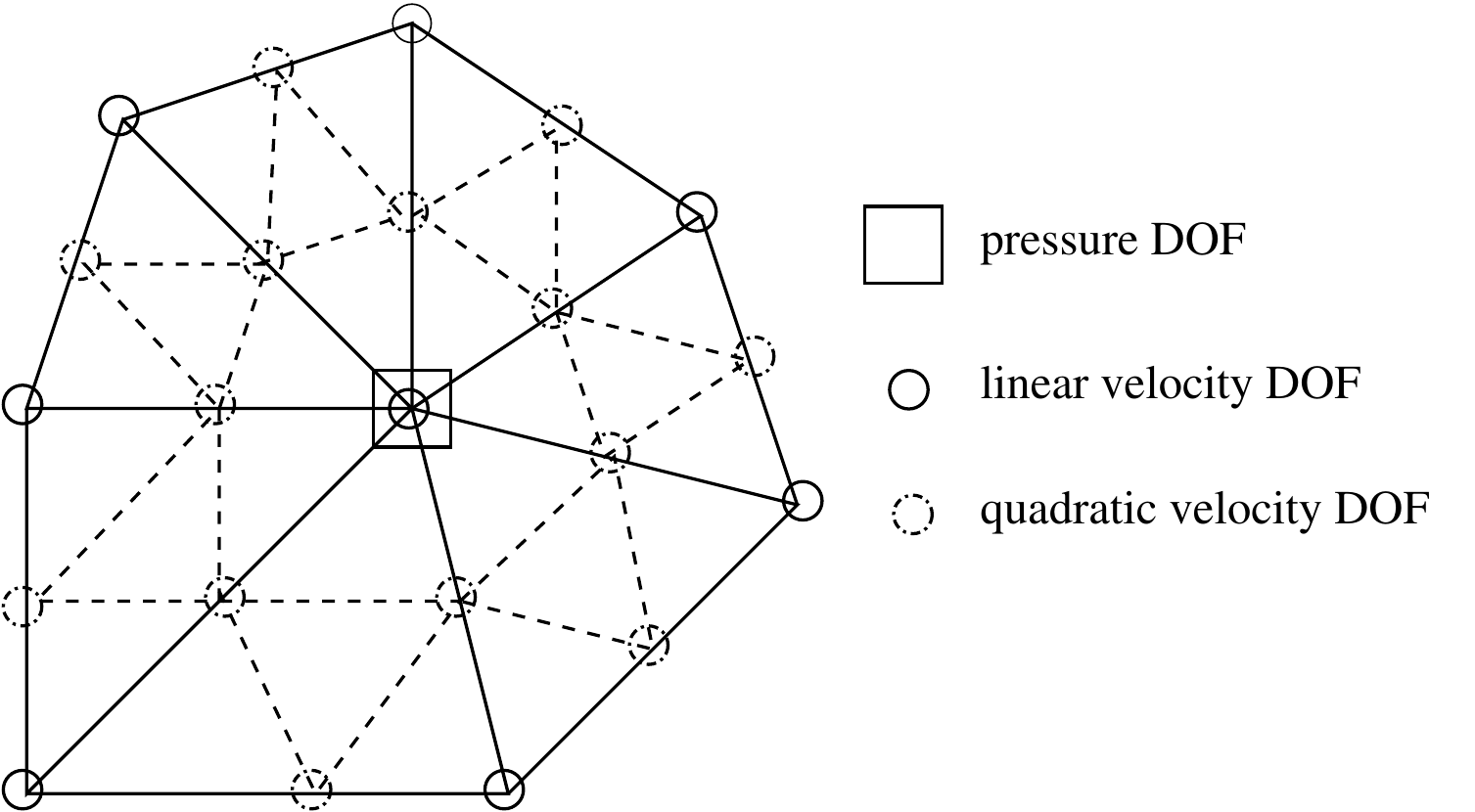}
  \caption{A typical local patch ${\mathcal P}_i$ 
    for the Taylor-Hood element contains one 
    pressure DOF, and connected linear and quadratic velocity DOF.}\label{fig:vankasm}
\end{figure}

The multiplicative Vanka smoother was introduced in \cite{Vanka86:00} for the 
fluid problem. We have recently developed an AMG method with 
this smoother for solving the nonlinear and nearly incompressible 
hyperelastic models in fluid-structure interaction 
simulation \cite{Langer201547}. To adapt this smoother for the Taylor-Hood element, we first 
construct the patches ${\mathcal P}_i$, $i=1, ... , n$. Each patch contains 
one pressure DOF, and the connected linear and quadratic velocity DOF 
indicated by the connectivity of matrices $B_{ll}$ and $B_{lq}$, respectively. A 
typical patch ${\mathcal P}_i$ is illustrated in Fig. \ref{fig:vankasm}. 

The local (correction) problem on ${\mathcal P}_i$ is extracted by a canonical 
projection of the global one (\ref{eq:mixedeqn}) to local one on ${\mathcal P}_i$: 
\begin{equation}\label{eq:locpatchsm}
  \left[ 
    \begin{array}{c}
      u_i^{k+1}\\
      p_i^{k+1}
    \end{array}
    \right]= \left[ 
  \begin{array}{c}
    u_i^{k}\\
    p_i^{k}
  \end{array}
  \right]+
    \omega
    \left[ 
    \begin{array}{cc}
      A_i & B_i^T\\
      B_i& -C_i
    \end{array}
    \right]^{-1}
  \left[ 
    \begin{array}{c}
      r_{u , i}^k\\
      r_{p, i}^k
    \end{array}
    \right]
\end{equation}
with $k$ representing the smoothing step and $\omega$ being a 
damping parameter. Here $[ (r_{u , i}^k)^T, (r_{p, i}^k)^T]^T$ denotes 
the residual updated in a multiplicative manner. As we observe from 
the numerical studies, this smoother shows the efficiency and 
robustness if it is used in the AMG preconditioner 
but not in the stand-alone AMG solver. 

\subsubsection{The Braess-Sarazin-type smoother}
This smoother was introduced in \cite{Braess97:00} and approximated in \cite{WZ00:00}, 
that has been applied to the fluid problem \cite{WM04:00, WM06:00}, 
the nearly incompressible elasticity problem \cite{YH11:00} and 
the fluid-structure interaction problem \cite{Yang20115367}. One smoothing 
step corresponds to a preconditioned Richardson method:
\begin{equation}\label{eq:bssm}
  \left[ 
    \begin{array}{c}
      u^{k+1}\\
      p^{k+1}
    \end{array}
    \right]= \left[ 
  \begin{array}{c}
    u^{k}\\
    p^{k}
  \end{array}
  \right]+\hat{K}^{-1}
  \left[ 
    \begin{array}{c}
      f - Au^k - B^T p^k \\
      g - Bu^k + Cp^k
    \end{array}
    \right]
\end{equation}
with preconditioner
\begin{equation}\label{eq:bssmpre}
  \hat{K}=
  \left[ 
    \begin{array}{cc}
      \hat{A} & B^T\\
      B & B\hat{A}^{-1}B^T - \hat{S}
    \end{array}
    \right].
\end{equation}
As in \cite{MW03:00}, we use $\hat{A}=2D$, where $D$ represents the diagonal of $A$. 
We use an AMG preconditioner $\hat{S}$ for the approximated 
Schur complement $C + B\hat{A}^{-1}B$. 

\subsubsection{The segregated Gauss-Seidel smoother}
This smoother was very recently introduced in \cite{Notay14:00} as 
a segregated Gauss-Seidel smoother based on a Uzawa-type iteration. 
Such a Uzawa method (see, e.g., \cite{ANU:298722}) 
can be reinterpreted as a preconditioned Richardson method: 
\begin{equation}\label{eq:uzsm}
  \left[ 
    \begin{array}{c}
      u^{k+1}\\
      p^{k+1}
    \end{array}
    \right]= \left[ 
  \begin{array}{c}
    u^{k}\\
    p^{k}
  \end{array}
  \right]+\hat{K}^{-1}
  \left[ 
   \begin{array}{c}
      f - Au^k - B^T p^k \\
     g - Bu^k + Cp^k
    \end{array}
    \right]
\end{equation}
with the preconditioner
\begin{equation}\label{eq:bssmpre}
  \hat{K}=
  \left[ 
    \begin{array}{cc}
      \hat{A} & \\
      B & -\omega^{-1} I 
    \end{array}
    \right],
\end{equation}
where $\omega$ is a properly chosen parameter, 
$\hat{A}$ is a (e.g., AMG) preconditioner for $A$. 
However, to get a multigrid smoother, 
this is relaxed by choosing some 
proper smoother $M_A$ for $A$ instead of a preconditioner $\hat{A}$; 
see \cite{Notay14:00}. 
The theoretical analysis requirement for $M_A$ has been specified therein. 
In our setting, we have chosen $M_A$ as a damped block Jacobi smoother with 
a damping parameter $0.5$. 
Compared to the Braess-Sarazin-type smoother, 
this smoother avoids the explicite construction of the approximate Schur complement. 

\section{Numerical results}\label{sec:numres}
\subsection{Meshes, boundary conditions and coarsening for 
the vector Laplacian and linear elasticity problems}
\begin{figure}[htbp]
  \centering
  \includegraphics[scale=0.28]{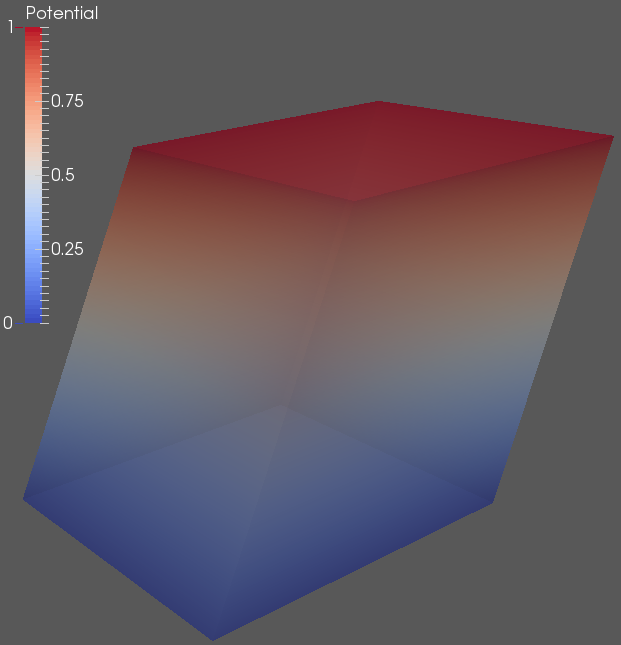}
  \includegraphics[scale=0.28]{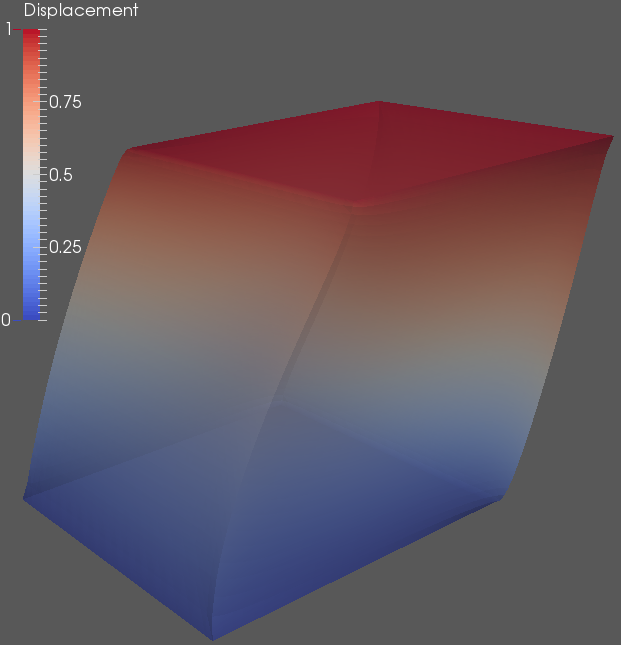}
  \caption{Numerical results of the vector Laplacian (left) and linear elasticity (right) 
    problem.}\label{fig:numres}
\end{figure}

We consider a unit cube $(0, 1)^3$ as the computational domain for 
the vector Laplacian and linear elasticity problems. The domain is subdivided 
into tetrahedra with four levels of mesh refinement $L_1-L_4$. 
The number of tetrahedron (\#Tet), nodes (\#Nodes) and midside 
nodes (\#Midside nodes), and the total number of DOF for elliptic 
( \#DOF (elliptic) ) and saddle point ((\#DOF (saddle point)) 
systems are shown in Table \ref{tab:meshinfo}.
\begin{table}[ht!]
   \centering
   \begin{tabular}{|c|c|c|c|c|}
     \hline
     Level  &  $L_1$ & $L_2$ & $L_3$ & $L_4$\\
     \hline
     \#Tet &  $64$ & $512$ & $4096$ & $32768$ \\
     \hline
     \#Nodes & $125$ &  $729$ &$4913$   &$35937$ \\
     \hline
     \#Midside nodes  & $604$ &$4184$ &$31024$   &$238688$ \\
     \hline
     \#DOF (elliptic)  & $2187$ &$14739$ &$107811$   &$823875$ \\
     \hline
     \#DOF (saddle point)  & $2312$ &$15468$ &$112724$   &$859812$ \\
     \hline
   \end{tabular}\caption{Number of tetrahedron (\#Tet) , nodes (\#Nodes), 
     Midside nodes (\#Midside nodes), and total number of DOF for 
     elliptic ( \#DOF (elliptic) ) and saddle point (\#DOF (saddle point)) systems 
     on four levels $L_1-L_4$.}\label{tab:meshinfo}
 \end{table} 
We fix the bottom of the domain, i.e., $u=[0, 0, 0]^T$ at $z=0$, prescribe 
a Dirichlet data on the top, i.e., $u=[0, 0, 1]^T$ at $z=1$, and 
use zero Neumann condition on the rest of the boundaries. For the linear 
elasticity problem, we set $\mu=1.15e+06$ and $\lambda=1.73e+06$. 
In Fig. \ref{fig:numres}, 
we plot the value of $\|u\|_{{\mathbb R}^3}$ of the 
numerical solution $u$ (indicated by the color) for the 
vector Laplacian (left) and the linear elasticity problem in pure displacement form 
(right), respectively. For visualization purpose, we plot the vector fields of potential 
and deformation, that is scaled by a factor of $0.1$. 

As a comparison, on each algebraically coarsening level, 
we show the number of linear and 
quadratic DOF (\# Linear DOF and \# Quadratic DOF, respectively)
in the new coarsening strategy, and 
the classical coarsening strategy (\# Non-separating DOF); see in Table \ref{tab:algdofs} the 
number of DOF on each coarsening level for the vector Laplacian 
and linear elasticity problem on the level $L_4$. It is obvious to see these two 
strategies lead to different graph connectivity in the coarsening procedure. 
 \begin{table}[ht!]
   \centering
   \begin{tabular}{|c|c|c|c|c|c|}
     \hline
     Coarsening Levels  &  $0$           &  $1$          & $2$        &  $3$      & $4$\\
     \hline
     \#Linear DOF           & $107811$   & $14739$   & $2187$   &  $375$  & $81$\\
     \hline
     \#Quadratic DOF      & $716064$   & $104544$ & $6366$   &  $369$  & $39$\\
     \hline 
     \#Non-separating DOF     & $823875$   & $88467$ & $2187$   &  $375$  & $24$\\
     \hline
   \end{tabular}\caption{Number of linear and quadratic DOF in the 
     new coarsening strategy, and the total DOF in the non-separating 
     coarsening strategy at level $L_4$, for the vector Laplacian and linear 
     elasticity problems.}\label{tab:algdofs}
  \end{table}

\subsection{Numerical performance for the vector Laplacian problem}
Before showing the AMG performance with the new 
coarsening strategy, we demonstrate the performance with a black-box type 
AMG \cite{FK98:00} in Table \ref{tab:nonsepamg} for the vector Laplacian problem, 
where non-separating coarsening strategy is used. 
For both AMG and AMG preconditioned CG methods, we use the relative 
residual error $\|f-Au_k\|_{l_2}/\|f-Au_0\|_{l_2}=1.0e-11$ in the $l_2-$norm as stopping 
criteria, where $k$ denotes the number of AMG iterations. 
In each iteration of the AMG solver, 
we use $2$ W-cycles and $1$ pre- and post-smoothing 
step. As observed, the AMG and the AMG preconditioned CG are not 
robust with respect to the mesh refinement, i.e., the iteration number 
increases with mesh refinement. 
 \begin{table}[ht!]
   \centering
   \begin{tabular}{|c|c|c|c|c|}
     \hline
      Levels  &  $L_1$           &  $L_2$          & $L_3$    &  $L_4$  \\
      \hline
      \#It AMG   & $90$   & $158$ & $>200$   &  $>200$\\
      \hline
      \#It PCG\_AMG      & $22$   & $25$   & $36$   &  $64$  \\
      \hline
   \end{tabular}\caption{Performance of black-box type AMG solver (\#It AMG) 
     and AMG preconditioned conjugate gradient solver (\#It PCG\_AMG) for the vector 
     Laplacian problem.}
   \label{tab:nonsepamg}
 \end{table}

Now we show the performance of 
the AMG and AMG preconditioned CG solvers using the new 
coarsening strategy. Note, that in the 
following numerical tests, we stop the iterations when the relative residual error 
in the $l_2-$norm is reduced by a factor $10^{11}$. 
We consider the 
Jacobi smoother with damping parameter $0.5$ and $1$ or $2$ 
pre- and post-smoothing steps (JA-1-1-0.5 or JA-2-2-0.5), 
and the Gauss-Seidel smoother with $1$ or $2$ 
pre- and post-smoothing steps (GS-1-1 or GS-2-2). In 
Table \ref{tab:amglapv}, we show the performance of the AMG solver 
for the vector Laplacian problem using V-cycle with different smoothers. 
In Table \ref{tab:amglapw}, we show the performance using W-cycle. 
In Table \ref{tab:amglapvpcg} and \ref{tab:amglapwpcg}, we show the 
performance of the AMG preconditioned CG using V- and W-cycles with 
different smoothers, respectively.

As observed, the iterations for each solver are independent of mesh refinement levels. 
The AMG solver using the Gauss-Seidel smoother shows better performance than 
the damped Jacobi smoother. By using the CG acceleration, we observe similar 
performance with two different smoothers. In addition, we observe, 
that the V- and W-cycles demonstrate almost the same performance. We also observe that 
the computational cost is proportional to the number of DOF. 
\begin{table}[ht!]
  \centering
  \begin{tabular}{|c|c|c|c|c|}
    \hline
    Level  &  $L_1$ & $L_2$ & $L_3$ & $L_4$\\
    \hline
    \#It ( JA-1-1-0.5 )  & $129$ & $125$ & $124$ & $128$  \\
    \hline
    \#It  ( JA-2-2-0.5 ) & $65$ & $65$ & $65$ & $66$  \\
    \hline
    \#It ( GS-1-1 )  & $43$ & $46$ & $47$ & $47$  \\
    \hline
    \#It  ( GS-2-2 ) & $23$ & $24$ & $24$ & $24$  \\
    \hline
  \end{tabular}
  \caption{Performance of the AMG solver for the vector Laplacian problem 
using V-cycle with different smoothers.}
  \label{tab:amglapv}
\end{table} 
\begin{table}[ht!]
  \centering
  \begin{tabular}{|c|c|c|c|c|}
    \hline
    Level  &  $L_1$ & $L_2$ & $L_3$ & $L_4$\\
    \hline
    \#It ( JA-1-1-0.5 )  & $128$ & $124$ & $121$ & $123$  \\
    \hline
    \#It  ( JA-2-2-0.5 ) & $65$ & $64$ & $63$ & $64$  \\
    \hline
    \#It ( GS-1-1 )  & $43$ & $46$ & $47$ & $47$  \\
    \hline
    \#It  ( GS-2-2 ) & $23$ & $24$ & $24$ & $24$  \\
    \hline
  \end{tabular}
  \caption{Performance of the AMG solver for the vector Laplacian problem 
using W-cycle with different smoothers.}
  \label{tab:amglapw}
\end{table} 
\begin{table}[ht!]
  \centering
  \begin{tabular}{|c|c|c|c|c|}
    \hline
    Level  &  $L_1$ & $L_2$ & $L_3$ & $L_4$\\
    \hline
    \#It ( JA-1-1-0.5 )  & $30$ & $30$ & $30$ & $30$  \\
    \hline
    \#It  ( JA-2-2-0.5 ) & $21$ & $22$ & $22$ & $22$  \\
    \hline
    \#It ( GS-1-1 )  & $26$ & $29$ & $29$ & $30$  \\
    \hline
    \#It  ( GS-2-2 ) & $17$ & $19$ & $19$ & $19$  \\
    \hline
  \end{tabular}
  \caption{Performance of the AMG preconditioned CG solver 
    for the vector Laplacian problem using V-cycle with different smoothers.}
  \label{tab:amglapvpcg}
\end{table} 
\begin{table}[ht!]
  \centering
  \begin{tabular}{|c|c|c|c|c|}
    \hline
    Level  &  $L_1$ & $L_2$ & $L_3$ & $L_4$\\
    \hline
    \#It ( JA-1-1-0.5 )  & $30$ & $30$ & $30$ & $29$  \\
    \hline
    \#It  ( JA-2-2-0.5 ) & $21$ & $21$ & $21$ & $21$  \\
    \hline
    \#It ( GS-1-1 )  & $26$ & $29$ & $29$ & $29$  \\
    \hline
    \#It  ( GS-2-2 ) & $17$ & $18$ & $18$ & $18$  \\
    \hline
  \end{tabular}
  \caption{Performance of the AMG preconditioned CG solver 
    for the vector Laplacian problem using W-cycle with different smoothers.}
  \label{tab:amglapwpcg}
\end{table} 

\subsection{Numerical performance for the linear elasticity problem 
  in pure displacement form}
We perform the same test for the linear elasticity problem. 
In Table \ref{tab:amgelasv}, we show the performance of the AMG solver 
for the linear elasticity problem using V-cycle with different smoothers. 
In Table \ref{tab:amgelasw}, we show the performance using W-cycle. 
In Table \ref{tab:amgelasvpcg} and \ref{tab:amgelaswpcg}, we show the 
performance of the AMG preconditioned CG using V- and W-cycles with 
different smoothers, respectively.

As observed, the damped Jacobi smoother does not work for this test problem. 
The AMG solver using the Gauss-Seidel smoother shows good performance. 
By using the CG acceleration, we observe improved performance. We observe, 
that the V- and W-cycles demonstrate almost the same performance.  
\begin{table}[ht!]
  \centering
  \begin{tabular}{|c|c|c|c|c|}
    \hline
    Level  &  $L_1$ & $L_2$ & $L_3$ & $L_4$\\
    \hline
    \#It ( JA-1-1-0.5 )  & $-$ & $-$ & $-$ & $-$  \\
    \hline
    \#It  ( JA-2-2-0.5 ) & $-$ & $-$ & $-$ & $-$  \\
    \hline
    \#It ( GS-1-1 )  & $80$ & $78$ & $76$ & $75$  \\
    \hline 
    \#It  ( GS-2-2 ) & $44$ & $40$ & $39$ & $39$  \\
    \hline
  \end{tabular}\caption{Performance of the AMG solver for the linear elasticity problem 
    in pure displacement form using V-cycle with different smoothers..}
  \label{tab:amgelasv}
\end{table} 
\begin{table}[ht!]
  \centering
  \begin{tabular}{|c|c|c|c|c|}
    \hline 
    Level  &  $L_1$ & $L_2$ & $L_3$ & $L_4$\\ \hline
    \#It ( JA-1-1-0.5 )  & $-$ & $-$ & $-$ & $-$  \\ \hline 
    \#It  ( JA-2-2-0.5 ) & $-$ & $-$ & $-$ & $-$  \\ \hline
    \#It ( GS-1-1 )  & $80$ & $78$ & $75$ & $73$  \\ \hline
    \#It  ( GS-2-2 ) & $44$ & $40$ & $39$ & $44$  \\
    \hline
  \end{tabular}\caption{Performance of the AMG solver for the linear elasticity problem 
    in pure displacement form using W-cycle with different smoothers..}
  \label{tab:amgelasw}
\end{table} 
\begin{table}[ht!]
  \centering
  \begin{tabular}{|c|c|c|c|c|}
    \hline
    Level  &  $L_1$ & $L_2$ & $L_3$ & $L_4$\\  \hline
    \#It ( JA-1-1-0.5 )  & $50$ & $81$ & $-$ & $-$  \\ \hline
    \#It  ( JA-2-2-0.5 ) & $32$ & $63$ & $-$ & $-$  \\ \hline
    \#It ( GS-1-1 )  & $40$ & $43$ & $43$ & $42$  \\ \hline
    \#It  ( GS-2-2 ) & $28$ & $27$ & $27$ & $27$  \\ 
    \hline
  \end{tabular}
  \caption{Performance of the AMG preconditioned CG solver 
    for the linear elasticity problem in pure displacement form 
    using V-cycle with different smoothers.}
  \label{tab:amgelasvpcg}
\end{table} 
\begin{table}[ht!]
  \centering
  \begin{tabular}{|c|c|c|c|c|}
    \hline 
    Level  &  $L_1$ & $L_2$ & $L_3$ & $L_4$\\ \hline
    \#It ( JA-1-1-0.5 )  & $47$ & $80$ & $-$ & $-$  \\ \hline
    \#It  ( JA-2-2-0.5 ) & $32$ & $62$ & $-$ & $-$  \\ \hline
    \#It ( GS-1-1 )  & $39$ & $44$ & $42$ & $41$  \\ \hline
    \#It  ( GS-2-2 ) & $28$ & $27$ & $27$ & $26$  \\
    \hline
  \end{tabular}
  \caption{Performance of the AMG preconditioned CG solver 
    for the linear elasticity problem in pure displacement form 
    using W-cycle with different smoothers.}
  \label{tab:amgelaswpcg}
\end{table} 

\subsection{Numerical performance for linear elasticity problem in mixed form}
\begin{figure}[htbp]
  \centering
  \includegraphics[scale=0.28]{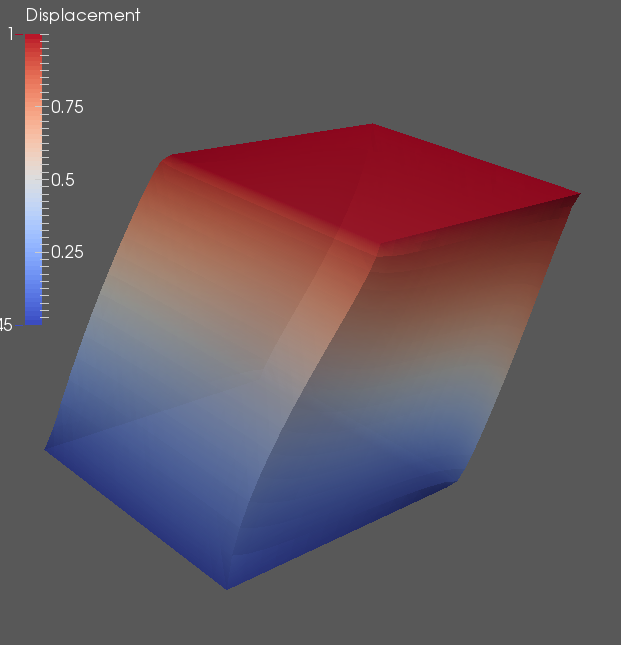}
  \includegraphics[scale=0.28]{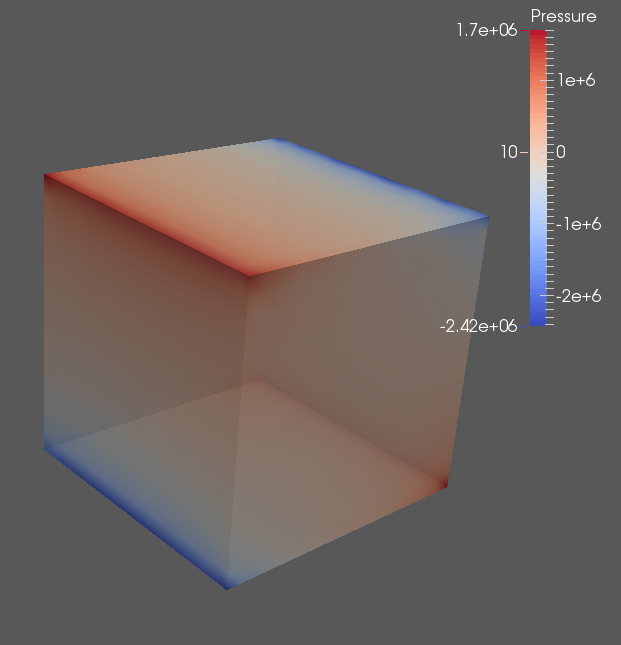}
  \caption{Numerical results of the displacement (left) and pressure (right) for 
    the linear elasticity problem in mixed form.}\label{fig:mixelasdispre}
\end{figure}

For the linear elasticity problem in mixed displacement-pressure form, we 
plot the simulation results of the displacement and pressure on the left and 
right plots of Fig. \ref{fig:mixelasdispre}, respectively. 

We set the 
relative residual error $1.0e-09$ in the corresponding norm as stopping criteria for solving the saddle point system 
(indefinite) with both the AMG solver and 
AMG preconditioned GMRES (see \cite{Saad86}) method. We will only 
consider the V-cyle for the remaining tests. 

The AMG solver with the Vanka smoother does not 
show the robustness and efficiency in this case. However combined with 
GMRES acceleration, the V-cycle preconditioner with such a smoother shows 
improved performance. We observe acceptable performance of 
one or two V-cycles (1 V-cycle or 2 V-cycle) preconditioned GMRES solver in 
Table \ref{tab:amgelasgmresvanka}. In each cycle, we use only one pre- and post-smoothing 
step.
\begin{table}[ht!]
  \centering
  \begin{tabular}{|c|c|c|c|c|}
    \hline 
    Level  &  $L_1$ & $L_2$ & $L_3$ & $L_4$\\ \hline
    \#It ( 1 V-cycle )  & $60$ & $45$ & $49$ & $69$  \\ \hline
    \#It  ( 2 V-cycles ) & $42$ & $30$ & $32$ & $45$  \\
    \hline
  \end{tabular}
  \caption{Performance of the V-cycle preconditioned GMRES solver 
    for the linear elasticity problem in mixed form using 
    one pre/poster Vanka smoother.}
  \label{tab:amgelasgmresvanka}
\end{table} 

The Braess-Sarazin smoother shows better performance. We observe 
the robustness with respect to the mesh size of the AMG solver and 
V-cycle preconditioned GMRES solver using such a smoother; see 
iteration numbers of the AMG solver using one or two 
Braess-Sarazin smoothing steps (Braess-Sarazin-1-1 or Braess-Sarazin-2-2) 
in Table \ref{tab:amgelasamgBS}, and one or two V-cycles (1 V-cycle or 2 V-cycle) 
preconditioned GMRES solver in Table \ref{tab:amgelasgmresBS}, respectively. 
\begin{table}[ht!]
  \centering
  \begin{tabular}{|c|c|c|c|c|}
    \hline 
    Level  &  $L_1$ & $L_2$ & $L_3$ & $L_4$\\ \hline
    \#It ( Braess-Sarazin-1-1 )  & $135$ & $133$ & $125$ & $117$  \\ \hline
    \#It  ( Braess-Sarazin-2-2 ) & $71$ & $70$ & $66$ & $62$  \\
    \hline
  \end{tabular}
  \caption{Performance of the AMG solver 
    for the linear elasticity problem in mixed form using 
    the V-cycle with the Braess-Sarazin smoother.}
  \label{tab:amgelasamgBS}
\end{table} 
\begin{table}[ht!]
  \centering
  \begin{tabular}{|c|c|c|c|c|}
    \hline 
    Level  &  $L_1$ & $L_2$ & $L_3$ & $L_4$\\ \hline
    \#It ( 1 V-cycle )  & $27$ & $27$ & $27$ & $27$  \\ \hline
    \#It  ( 2 V-cycles ) & $18$ & $19$ & $19$ & $19$  \\
    \hline
  \end{tabular}
  \caption{Performance of the V-cycle preconditioned GMRES solver 
    for the linear elasticity problem in mixed form using one pre/post Braess-Sarazin smoother.}
  \label{tab:amgelasgmresBS}
\end{table} 

Using the segregated Gauss-Seidel smoother (sGS), we observe good performance. 
For all tests, we use $\omega=0.125$. 
The robustness with respect to the mesh size of the AMG solver can be observed; see 
iteration numbers of the AMG solver using one or two 
segregated Gauss-Seidel smoothing steps (sGS-1-1 or sGS-2-2) 
in Table \ref{tab:amgelasamguzawanotay}. The efficiency is further improved when 
combined with the Krylov subspace acceleration; see iteration numbers of 
one or two V-cycles (1 V-cycle or 2 V-cycle) 
preconditioned GMRES solver in Table \ref{tab:amgelasgmresuzawanotay}.
\begin{table}[ht!]
  \centering
  \begin{tabular}{|c|c|c|c|c|}
    \hline 
    Level  &  $L_1$ & $L_2$ & $L_3$ & $L_4$\\ \hline
    \#It ( sGS-1-1 )  & $107$ & $104$ & $99$ & $93$  \\ \hline
    \#It  ( sGS-2-2 ) & $54$ & $53$ & $53$ & $59$  \\
    \hline
  \end{tabular}
  \caption{Performance of the AMG solver 
    for the linear elasticity problem in mixed form using 
    the V-cycle with the segregated Gauss-Seidel smoother.}
  \label{tab:amgelasamguzawanotay}
\end{table} 
\begin{table}[ht!]
  \centering
  \begin{tabular}{|c|c|c|c|c|}
    \hline 
    Level  &  $L_1$ & $L_2$ & $L_3$ & $L_4$\\ \hline
    \#It ( 1 V-cycle )  & $14$ & $18$ & $19$ & $21$  \\ \hline
    \#It  ( 2 V-cycles ) & $12$ & $15$ & $15$ & $15$  \\
    \hline
  \end{tabular}
  \caption{Performance of the V-cycle preconditioned GMRES solver 
    for the linear elasticity problem in mixed form using one pre/post segregated 
    Gauss-Seidel smoother.}
  \label{tab:amgelasgmresuzawanotay}
\end{table} 

\subsection{Numerical performance for the Stokes problem}
\begin{figure}[htbp]
  \centering
  \includegraphics[scale=0.26]{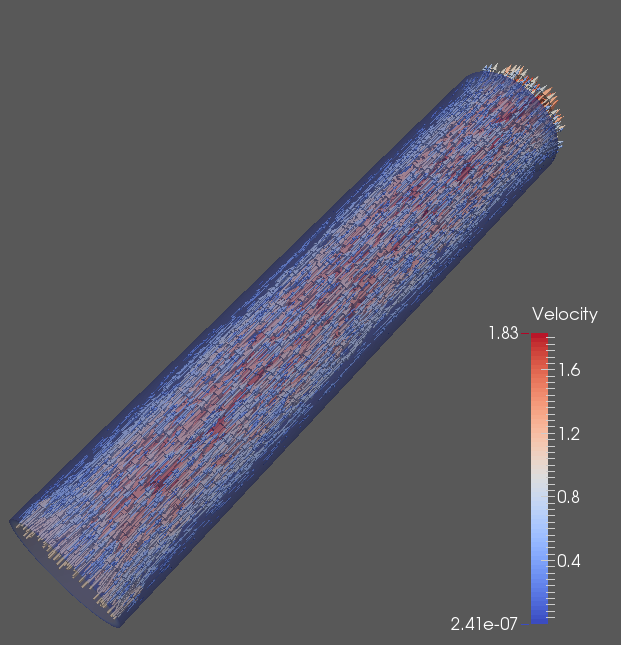}
  \includegraphics[scale=0.26]{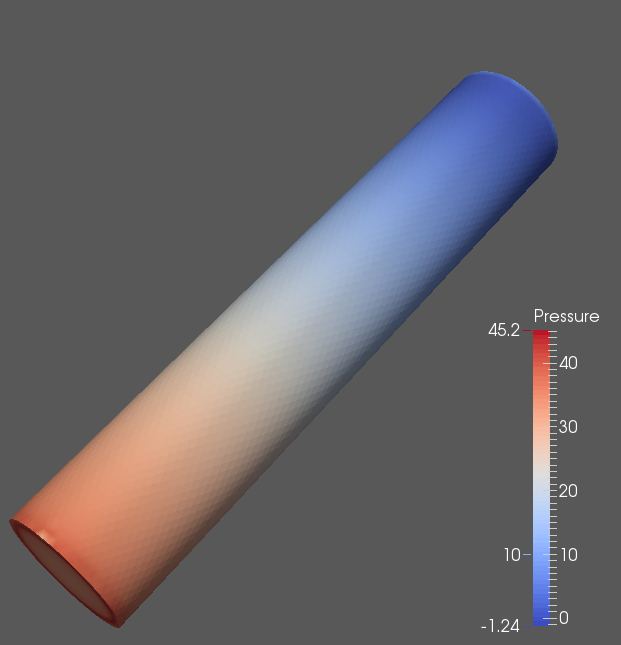}
  \caption{Numerical results of the Stokes velocity (left) and pressure (right).}
  \label{fig:stokesnumres}
\end{figure}

The computational domain for the Stokes problem is prescribed by an 
inside of a cylinder, that has radius of $1$ with center point $(0, 0, 0 )$ 
on the inflow boundary (where $u=(1.0, 0, 0)$), 
and center point $(10, 0, 0)$ on the outflow 
boundary (where $(2\mu\varepsilon(u)-pI)n=(0, 0, 0)$). On the 
rest of the boundaries $u=(0, 0, 0)$. For all tests, we set $\mu=0.5$. 
Four levels ($L_1-L_4$) of tetrahedral meshes are generated; see mesh information for 
each level in Table \ref{tab:fluidmeshinfo}: The number of tetrahedron (\#Tet), 
nodes (\#Nodes) and midside nodes (\#Midside nodes), and the total number 
of DOF ( \#DOF ) for the saddle point system. As an illustration to show the coarsening 
strategy, in Table \ref{tab:algdofsfluid}, we show the number of linear and 
quadratic velocity DOF ( \# Linear velocity DOF and  \# Quadratic velocity DOF) , 
the linear pressure DOF (\# Linear pressure DOF) on each coarsening level for the level 
$L_4$. The velocity and pressure of the simulation results are shown on the left and right 
plots in Fig. \ref{fig:stokesnumres}, respectively.  
\begin{table}[ht!]
   \centering
   \begin{tabular}{|c|c|c|c|c|}
     \hline
     Level  &  $L_1$ & $L_2$ & $L_3$ & $L_4$\\
     \hline
     \#Tet &  $895$ & $7160$ & $57280$ & $458240$ \\
     \hline
     \#Nodes & $351$ &  $1922$ &$12307$   &$87109$ \\
     \hline
     \#Midside nodes  & $1571$ &$10385$ &$74802$   &$566212$ \\
     \hline
     \#DOF & $6117$ &$38843$ &$273634$   &$2047072$ \\
     \hline
   \end{tabular}\caption{Fluid mesh information: Number of tetrahedron (\#Tet) , 
     nodes (\#Nodes) 
     and Midside nodes (\#Midside nodes), and the total number of DOF (\#DOF) 
     on four levels $L_1-L_4$.}\label{tab:fluidmeshinfo}
 \end{table} 
 \begin{table}[ht!]
   \centering
   \begin{tabular}{|c|c|c|c|c|c|}
     \hline
     Coarsening Levels  &  $0$           &  $1$          & $2$        &  $3$      & $4$\\
     \hline
     \#Linear velocity DOF           & $261327$   & $36921$   & $5766$   &  $1053$  & $267$\\
     \hline
     \#Quadratic velocity DOF      & $1698636$   & $213492$ & $14142$   &  $1053$  & $93$\\
     \hline 
     \#Linear pressure DOF           & $87109$   & $12307$   & $1922$   &  $351$  & $89$\\
     \hline
   \end{tabular}\caption{Number of linear and quadratic velocity, and linear 
     pressure DOF in the new coarsening strategy at level $L_4$ 
     for the Stokes problem.}\label{tab:algdofsfluid}
  \end{table}

For this example, the AMG solver with the 
Vanka, Braess-Sarazin and segregated Gauss-Seidel smoothers shows 
poor performance, that is very large smoothing steps are required in order to get multigrid convergence rate. 
However, this will lead to very expensive computational cost. 
In addition, we observe unsatisfactory performance of the AMG preconditioned 
Krylov subspace method using the V-cycle with the Vanka and 
segregated Gauss-Seidel smoothers. Therefore, we only report the performance of 
the AMG preconditioned GMRES solver using the Braess-Sarazin 
smoother, that is shown in Table 
\ref{tab:amgstokesgmresBS}. We set relative residual error $1.0e-09$ in the corresponding 
norm as stopping criteria of the AMG preconditioned GMRES solver. 
It is easy to see, with the Krylov subspace 
acceleration, the performance is greatly improved, using one or two V-cycle 
(1 V-cycle or 2 V-cycle) preconditioner with one pre- and post-smoothing steps. 
\begin{table}[ht!]
  \centering
  \begin{tabular}{|c|c|c|c|c|}
    \hline 
    Level  &  $L_1$ & $L_2$ & $L_3$ & $L_4$\\ \hline
    \#It ( 1 V-cycle )  & $25$ & $42$ & $38$ & $39$  \\ \hline
    \#It  ( 2 V-cycles ) & $16$ & $17$ & $18$ & $21$  \\
    \hline
  \end{tabular}
  \caption{Performance of the V-cycle preconditioned GMRES solver 
    for the Stokes problem using one pre/post Braess-Sarazin smoother.}
  \label{tab:amgstokesgmresBS}
\end{table} 

\section{Conclusions}\label{sec:con}
In this work, we have developed an AMG method used as a stand-alone 
solver or preconditioner in the Krylov subspace methods for solving the 
finite element equations of the vector Laplacian problem, linear elasticity problem in 
pure displacement and mixed displacement-pressure form, and the Stokes 
problem in mixed velocity-pressure form in 3D. 
We have developed a new strategy to construct the 
hierarchy of the AMG coarsening system using the 
hierarchical quadratic basis functions. The numerical studies have demonstrated 
the good performance of the AMG solvers or the AMG preconditioned 
Krylov subspace methods for the elliptic and saddle point systems, respectively. 
In particular, the AMG preconditioned Krylov subspace methods show much 
better robustness and efficiency for solving both systems 
compared with the AMG stand-alone solvers. 
From this point of view, the AMG method developed in this work can be used 
as a robust and efficient solver or preconditioner 
for the SPD system and the saddle point system with compressible materials, and as 
a robust and efficient preconditioner for the saddle point 
system with incompressible materials. It is also possible to extend this AMG method 
for high-order hierarchical finite element basis functions. 

\section*{Acknowledgement} 
The author would like to thank Prof. Ulrich Langer for his encouragement 
and many enlightened discussions on this work. 

\bibliography{amgp2p1}

\newcommand{\noopsort}[1]{} \newcommand{\printfirst}[2]{#1}
  \newcommand{\singleletter}[1]{#1} \newcommand{\switchargs}[2]{#2#1}
\begin{thebibliography}{10}

\bibitem{AF03:00}
R.A. Adams and J.J.F. Fournier.
\newblock {\em Sobolev Spaces}.
\newblock Academic Press, Amsterdam, Boston, 2003.

\bibitem{ANU:298722}
Michele Benzi, G.H. Golub, and J.~Liesen.
\newblock Numerical solution of saddle point problems.
\newblock {\em Acta Numerica}, 14:1--137, 5 2005.

\bibitem{DB95}
D.~Braess.
\newblock Towards algebraic multigrid for elliptic problems of second order.
\newblock {\em Computing}, 55(4):379--393, 1995.

\bibitem{DB07:00}
D.~Braess.
\newblock {\em Finite Elements - Theory, Fast Solvers, and Applications in
  Solid Mechanics}.
\newblock Cambridge University Press, Cambridge, New York, 2007.

\bibitem{Braess97:00}
D.~Braess and R.~Sarazin.
\newblock {An efficient smoother for the Stokes problem}.
\newblock {\em Appl. Numer. Math.}, 23(1):3--19, 1997.

\bibitem{BF91:00}
F.~Brezzi and M.~Fortin.
\newblock {\em {Mixed and Hybrid Finite Element Methods}}.
\newblock Springer, New York, 1991.

\bibitem{Notay14:00}
F.J. Gaspar, Y.~Notay, C.W. Oosterlee, and C.~Rodrigo.
\newblock {A simple and efficient segregated smoother for the discrete Stokes
  equations}.
\newblock {\em SIAM J. Sci. Comput.}, 36(3):A1187--A1206, 2014.

\bibitem{UH:02}
G.~Haase and U.~Langer.
\newblock {\em Modern Methods in Scientific Computing and Applications},
  volume~75 of {\em NATO Science Series II. Mathematics, Physics and
  Chemistry}, chapter Multigrid Methods: From Geometrical to Algebraic
  Versions, pages 103--154.
\newblock Kluwer Academic Press, Dordrecht, 2002.

\bibitem{HB03}
W.~Hackbusch.
\newblock {\em Multi-Grid Methods and Applications}.
\newblock Springer, Heidelberg, 2003.

\bibitem{JA08:00}
A.~Janka.
\newblock Smoothed aggregation multigrid for a {S}tokes problem.
\newblock {\em Comput. Visual. Sci.}, 11(3):169--180, 2008.

\bibitem{FK98:00}
F.~Kickinger.
\newblock Algebraic multigrid for discrete elliptic second-order problems.
\newblock In {\em Multigrid Methods V. Proceedings of the 5th European
  Multigrid conference (ed.~by W. Hackbush), Lecture Notes in Computational
  Sciences and Engineering, vol.~3}, pages 157--172. Springer, 1998.

\bibitem{Langer2005406}
U.~Langer and D.~Pusch.
\newblock Data-sparse algebraic multigrid methods for large scale boundary
  element equations.
\newblock {\em Appl. Numer. Math.}, 54(3–4):406--424, 2005.

\bibitem{NME:NME839}
U.~Langer, D.~Pusch, and S.~Reitzinger.
\newblock Efficient preconditioners for boundary element matrices based on
  grey-box algebraic multigrid methods.
\newblock {\em Int J Numer Meth Engng}, 58(13):1937--1953, 2003.

\bibitem{Langer201547}
U.~Langer and H.~Yang.
\newblock Partitioned solution algorithms for fluid-structure interaction
  problems with hyperelastic models.
\newblock {\em J. Comput. Appl. Math.}, 276(0):47--61, 2015.

\bibitem{BM:13}
B.~Metsch.
\newblock {\em {Algebraic Multigrid (AMG) for Saddle Point Systems}}.
\newblock PhD thesis, Rheinischen Friedrich-Wihelms-Universit\"{a}t Bonn, 2013.

\bibitem{NOT14}
A.~Napov and Y.~Notay.
\newblock Algebraic multigrid for moderate order finite elements.
\newblock {\em SIAM J Sci Comput}, 2014.
\newblock to appear.

\bibitem{SR01:00}
S.~Reitzinger.
\newblock {\em {Algebraic Multigrid Methods for Large Scale Finite Element
  Methods}}.
\newblock PhD thesis, Johannes Kepler University Linz, 2001.

\bibitem{JWRKS87}
J.~W. Ruge and K.~St{\"u}ben.
\newblock Algebraic multigrid.
\newblock In S.F. McCormick, editor, {\em Multigrid Methods}, volume~3 of {\em
  Frontiers in Applied Mathematics}, pages 73--130. SIAM, Philadelphia, PA,
  1987.

\bibitem{Saad86}
Y.~Saad and Martin~H. Schultz.
\newblock {GMRES}: {A} generalized minimal residual algorithm for solving
  nonsymmetric linear systems.
\newblock {\em SIAM J. Sci. Stat. Comput.}, 7(3):856--869, 1986.

\bibitem{XUJ06}
S.~Shu, D.~Sun, and J.~Xu.
\newblock An algebraic multigrid method for higher-order finite element
  discretizations.
\newblock {\em Computing}, 77(4):347--377, 2006.

\bibitem{KS01:00}
K.~St\"{u}ben.
\newblock {\em Multigrid}, chapter Appendix A: An introduction to algebraic
  multigrid, pages 413--533.
\newblock Academic Press, 2001.

\bibitem{Stuben2001281}
K.~St\"{u}ben.
\newblock A review of algebraic multigrid.
\newblock {\em J. Comput. Appl. Math.}, 128(1–2):281--309, 2001.

\bibitem{Vanka86:00}
S.P. Vanka.
\newblock {Block-implicit multigrid solution of Navier-Stokes equations in
  primitive variables}.
\newblock {\em J. Comput. Phys.}, 65(1):138--158, 1986.

\bibitem{MW03:00}
M.~Wabro.
\newblock {\em {Algebraic Multigrid Methods for the Numerical Solution of the
  Incompressible Navier-Stokes Equations}}.
\newblock PhD thesis, Johannes Kepler University Linz, 2003.

\bibitem{WM04:00}
M.~Wabro.
\newblock Coupled algebraic multigrid methods for the {O}seen problem.
\newblock {\em Comput Visual Sci}, 7(3-4):141--151, 2004.

\bibitem{WM06:00}
M.~Wabro.
\newblock {AMGe}---coarsening strategies and application to the {O}seen
  equations.
\newblock {\em SIAM J Sci Comput}, 27(6):2077--2097, 2006.

\bibitem{TW15}
T.~Wiesner.
\newblock {\em {Flexible Aggregration-based Algebraic Multigrid Method for
  Contact and Flow Problems}}.
\newblock PhD thesis, Technischen Universit\"{a}t M\"{u}nchen, 2015.

\bibitem{YH11:00}
H.~Yang.
\newblock Partitioned solvers for the fluid-structure interaction problems with
  a nearly incompressible elasticity model.
\newblock {\em Comput. Visual. Sci.}, 14(5):227--247, 2011.

\bibitem{Yang20115367}
H.~Yang and W.~Zulehner.
\newblock Numerical simulation of fluid-structure interaction problems on
  hybrid meshes with algebraic multigrid methods.
\newblock {\em J. Comput. Appl. Math.}, 235(18):5367--5379, 2011.

\bibitem{WZ00:00}
W.~Zulehner.
\newblock A class of smoothers for saddle point problems.
\newblock {\em Computing}, 65(3):227--246, 2000.

\end{thebibliography}
\bibliographystyle{plain}

\end{document}